\input amstex
\documentstyle{amsppt}
\magnification 1200
\NoBlackBoxes
\pageheight{48pc}
\topmatter
\title{ON THE BERNSTEIN-GEL'FAND-GEL'FAND CORRESPONDENCE AND A RESULT OF 
EISENBUD, FL\O YSTAD, AND SCHREYER}\endtitle
\rightheadtext{ON THE BGG CORRESPONDENCE}
\leftheadtext{IUSTIN COAND\u A}
\author{IUSTIN COAND\u A}\endauthor
\address{INSTITUTE OF MATHEMATICS OF THE ROMANIAN ACADEMY, P.O.BOX 1-764, 
RO-70700 BUCHAREST, ROMANIA}
\endaddress
\email{Iustin.Coanda\@imar.ro}\endemail
\abstract{We show that a combination between a remark of I.N. Bernstein, 
I.M. Gel'fand, and S.I. Gel'fand [2] and the idea, systematically investigated 
by D. Eisenbud, G. Fl\o ystad, and F.-O. Schreyer [3], of taking Tate 
resolutions over exterior algebras leads to quick proofs of the main results 
of [2] and [3] (theorems 7 and 10 below). This combination is expressed by 
lemma 6 from the text, a result which we prove directly using the 
cohomology of invertible sheaves on a projective space.}
\endabstract
\endtopmatter
\document

Since the above abstract may serve as an introduction as well, we begin by 
recalling (in (0)-(4)) some definitions and facts. We use the Chapter I of [5] 
as our main reference for homological algebra (except that we denote mapping 
cones by ``Con''). 

\subhead{0. Definition }\endsubhead  
Let $k$ be a field, $V$ an $(n+1)$-dimensional 
$k$-vector space, $e_0,...,e_n$ a $k$-basis of $V$ and $X_0,...,X_n$ the 
dual basis of $V^*$. Let $\Lambda =\wedge (V)$ be the exterior algebra of
$V$. $\Lambda$ is a (positively) graded $k$-algebra : 
$\Lambda ={\Lambda}_0\oplus ...\oplus {\Lambda}_{n+1}$ with 
${\Lambda}_i={\wedge}^i(V)$. Let ${\Lambda}_+:={\Lambda}_1\oplus ...\oplus 
{\Lambda}_{n+1}$ and $\underline k:={\Lambda}/{\Lambda}_+$. 
We denote by $\Lambda$-mod the category of finitely generated, graded, 
right $\Lambda$-modules (with morphisms of degree 0). 

Let $\Bbb P=\Bbb P(V)$ be the projective space of 1-dimensional $k$-vector 
subspaces of $V$ (such that $H^0{\Cal O}_{\Bbb P}(1)=V^*$). If 
$N\in \text{Ob}(\Lambda \text{-mod})$ one defines a bounded complex 
$\text{L}(N)$ of coherent sheaves on $\Bbb{P}(V)$ by  
$\text{L}(N)^p:={\Cal O}_{\Bbb P}(p){\otimes}_kN_p$ and 
$d_{{\text L}(N)}:={\mathop{\sum}\limits _{i=0}^n}(X_i\cdot -)\otimes 
(-\cdot e_i)$. In this way one obtains {\it the BGG functor }\rm 
$\text{L} : \Lambda \text{-mod}\rightarrow \text{C}^b(\text{Coh}\Bbb{P}(V))$. 
It can be extended to a functor $\text{L} : \text{C}(\Lambda \text{-mod})
\rightarrow \text{C}(\text{Qcoh}\Bbb{P}(V))$ as it follows : if 
$K^{\bullet}$ is a complex in $\Lambda$-mod one considers the double complex 
$X^{\bullet \bullet}$ in Coh$\Bbb{P}(V)$ with 
$X^{p,\bullet}:=\text{L}(K^p)$ and with 
${d'}^p : X^{p,\bullet}\rightarrow X^{p+1,\bullet}$ equal to  
$\text{L}(d_K^p)$ and one takes 
$\text{L}(K^{\bullet}):=\text{s}(X^{\bullet \bullet})$ the simple complex 
associated to $X^{\bullet \bullet}$. 

The (extended) functor L is exact, commutes with the translation functor T and 
with mapping cones and maps morphisms homotopically equivalent to 0 to 
morphisms with the same property (see [3] remark after (2.5) for a nice 
argument) hence it induces a functor 
L : $\text{K}(\Lambda \text{-mod})\rightarrow \text{K}(\text{Qcoh}\Bbb{P}(V))$.
L also maps quasi-isomorphisms in $\text{K}^+(\Lambda \text{-mod})$ to 
quasi-isomorphisms in $\text{K}(\text{Qcoh}\Bbb{P}(V))$, hence it induces a 
functor L : $\text{D}^+(\Lambda \text{-mod})\rightarrow \text{D}(\text{Qcoh}
\Bbb{P}(V))$. 

We shall often use the following shorter notations : 
$\text{K}(\Lambda ):=\text{K}(\Lambda \text{-mod})$, 
$\text{D}(\Lambda ):=\text{D}(\Lambda \text{-mod})$, 
$\text{D}(\Bbb P):=\text{D}(\text{Qcoh}\Bbb{P}(V))$ and 
$\text{D}^b(\Bbb P):=\text{D}^b(\text{Coh}\Bbb{P}(V))$. 

\subhead{1. Definition }\endsubhead
(i) If $N\in \text{Ob}(\Lambda \text{-mod})$ and $a\in \Bbb{Z}$ one defines 
a new object $N(a)$ of $\Lambda $-mod by : $N(a)_p:=N_{a+p}$ and 
$(y\cdot v)_{N(a)}:=(-1)^a(y\cdot v)_N,\  \forall y\in N,\  \forall v\in V$. 
With 
this convention, if $\omega \in {\Lambda}_b$ then $(-\cdot \omega)_N$ defines 
a morphism in $\Lambda $-mod : $N(a)\rightarrow N(a+b)$. If 
$u : N'\rightarrow N$ is a morphism then $u(a) : N'(a)\rightarrow N(a)$ is 
just $u$ if one forgets the gradings.

One has : $\text{L}(N(a))=\text{T}^a\text{L}(N)(-a)$. If $K^{\bullet}$ is a 
complex in $\Lambda$-mod, let $K^{\bullet}((a))$ be the complex which 
coincides with $K^{\bullet}(a)$ term by term but with 
$d_{K((a))}:=(-1)^ad_{K(a)}$. Then 
$\text{L}(K^{\bullet}((a)))=\text{T}^a\text{L}(K^{\bullet})(-a)$ and if one 
applies L to the isomorphism 
$((-1)^{ap}\cdot \text{id}_{K^p(a)})_{p\in \Bbb{Z}} : 
K^{\bullet}(a)\buildrel \sim \over \rightarrow K^{\bullet}((a))$ one gets 
a functorial isomorphism 
$\text{L}(K^{\bullet}(a))\simeq \text{T}^a\text{L}(K^{\bullet})(-a)$.

\vskip 0.1in

(ii) If $N\in \text{Ob}(\Lambda \text{-mod})$ let $N\spcheck $ denote the 
graded $k$-vector space $\text{Hom}_k(N,k)$ endowed with the following 
right $\Lambda $-module structure : for $v\in V$, the multiplication 
$(-\cdot v)_{N\spcheck } : (N\spcheck )_p\rightarrow (N\spcheck )_{p+1}$ is, 
by definition, $(-1)^{p+1}\cdot $the dual of the multiplication 
$(-\cdot v)_N : N_{-p-1}\rightarrow N_{-p}$. With this definition, 
$\text{L}(N\spcheck )={\Cal Hom}_{{\Cal O}_\Bbb P}^{\bullet}(\text{L}(N),
{\Cal O}_{\Bbb P})$. 

The canonical isomorphism of $k$-vector spaces 
$\mu : N\rightarrow (N\spcheck ) \spcheck $ is not a morphism in 
$\Lambda $-mod : $\mu (y\cdot v)=-\mu (y)\cdot v,\  \forall y\in N, 
\  \forall v\in V$. However, ${\mu}':=((-1)^p{\mu}_p)_{p\in \Bbb{Z}}$ defines 
an 
isomorphism in $\Lambda $-mod : $N\buildrel \sim \over \rightarrow 
(N\spcheck ) \spcheck$.

\vskip 0.1in

(iii) Of a particular importance is the object $\Lambda \spcheck $ of 
$\Lambda $-mod. One has $(\Lambda \spcheck)_{-p}={\wedge}^pV^*,\  \forall 
p\in \Bbb{Z}$ and, for $v\in V$, the multiplication 
$(-\cdot v)_{\Lambda \spcheck}: (\Lambda \spcheck)_{-p}\rightarrow 
(\Lambda \spcheck)_{-p+1}$ is the contraction by $v$ : 
$(f_1\wedge ...\wedge f_p\cdot v)_{\Lambda \spcheck}=
\mathop{\sum}\limits _{i=1}^p(-1)^{i-1}f_i(v)\cdot 
f_1\wedge ...\wedge {\widehat{f_i}}\wedge ...\wedge f_p$ for 
$f_1,...,f_p\in V^*$. It follows that $\text{L}(\Lambda \spcheck)$ is the 
tautological Koszul complex on $\Bbb{P}(V)$ : 

$$ 
0\rightarrow {\Cal O}_{\Bbb P}(-n-1){\otimes}_k{\wedge}^{n+1}V^*
\rightarrow ...\rightarrow 
{\Cal O}_{\Bbb P}(-1){\otimes}_kV^*\rightarrow {\Cal O}_{\Bbb P}@>>>0.
$$

(iv) If $N\in \text{Ob}(\Lambda \text{-mod})$, soc($N$) consists of the 
elements of $N$ annihilated by ${\Lambda}_+$. In particular, 
soc($\Lambda $)=${\Lambda}_{n+1}$ and 
soc($\Lambda \spcheck$)=$(\Lambda \spcheck)_0$. 

\subhead{2. Remark }\endsubhead 
(i) Let $\Cal A$ be an abelian category. Consider a short exact sequence :   
$$ 
0\rightarrow X^{\bullet}\overset u\to\longrightarrow Y^{\bullet}
\overset v\to\longrightarrow Z^{\bullet}\rightarrow 0
$$
in the category $\text{C}(\Cal A)$ 
of complexes in $\Cal A$. Let 
$w : Z^{\bullet}\rightarrow \text{T}X^{\bullet}$ be the morphism in the 
derived category $\text{D}(\Cal A)$ defined by the diagram : 
$$\CD 
Z^{\bullet}@<(0,v)<\text{qis}<\text{Con}(u)@>(\text{id}_{\text{T}X},0)>>
\text{T}X^{\bullet}
\endCD 
$$
(recall that $\text{Con}(u)=\text{T}X^{\bullet}\oplus Y^{\bullet}$ term 
by term, not as complexes). Then 
$(X^{\bullet},Y^{\bullet},Z^{\bullet}$ $,u,v,w)$ is a 
$distinguished$ $triangle$ in $\text{D}(\Cal A)$ hence 
$(Y^{\bullet},Z^{\bullet},\text{T}X^{\bullet},v,w,-\text{T}u)$ and 
$(\text{T}^{-1}Z^{\bullet},X^{\bullet},Y^{\bullet},$ $-\text{T}^{-1}w,u,v)$ 
are distiguished triangles too. One gets a ``long'' complex in 
$\text{D}(\Cal A)$ : 
$$\CD 
...\longrightarrow \text{T}^{-1}Z^{\bullet}@>-\text{T}^{-1}w>>X^{\bullet}
\overset u\to\longrightarrow Y^{\bullet}\overset v\to\longrightarrow 
Z^{\bullet}\overset w\to\longrightarrow \text{T}X^{\bullet}
@>-\text{T}u>>\text{T}Y^{\bullet}\longrightarrow ...
\endCD 
$$
and for every $W^{\bullet}\in \text{ObC}(\Cal A)$ if one applies 
$\text{Hom}_{\text{D}(\Cal A)}(W^{\bullet},-)$ or 
$\text{Hom}_{\text{D}(\Cal A)}(-,W^{\bullet})$ to this long complex one 
gets long exact sequences in the category $\Cal Ab$ of abelian groups. 

\vskip 0.1in

(ii) Assume that the short exact sequence of complexes from (i) is 
$semi$-$split$, i.e., that 
$$\CD 
0\longrightarrow X^p@>u^p>>Y^p@>v^p>>Z^p\longrightarrow 0 
\endCD 
$$ 
is split-exact $\forall p\in \Bbb Z$, i.e., there exist morphisms 
$s^p : Z^p\rightarrow Y^p$ and $t^p : Y^p\rightarrow X^p$ such that 
$t^p\circ u^p=\text{id}_{X^p}$, $v^p\circ s^p=\text{id}_{Z^p}$ and 
$\  u^p\circ t^p+s^p\circ v^p=\text{id}_{Y^p}$ (hence 
$\  t^p\circ s^p=0$). Then 
$\delta :=(t^{p+1}\circ d_Y^p\circ s^p)_{p\in \Bbb Z}$ is a morphism of 
complexes (i.e., in $\text{C}(\Cal A)$) : 
$Z^{\bullet}\rightarrow \text{T}X^{\bullet}$ and $w=-\delta $ in 
$\text{D}(\Cal A)$ (in fact, 
$(^t(-{\delta}^p,s^p))_{p\in \Bbb Z} : Z^{\bullet}\rightarrow 
\text{Con}(u)$ is an inverse of $(0,v)$ in 
$\text{K}(\Cal A)$). 

\vskip 0.1in

(iii) Assume that $X^{\bullet}, Y^{\bullet}, Z^{\bullet}\in 
\text{ObC}^+(\Cal A)$ and consider a short exact sequence as in (i). 
Let $I^{\bullet}\in \text{ObC}(\Cal A)$ be a complex consisting of 
$injective$ objects of $\Cal A$. Then the functor 
$\text{Hom}_{\text{K}(\Cal A)}(-,I^{\bullet})$ maps quasi-isomorphisms 
in $\text{K}^+(\Cal A)$ to isomorphisms in $\Cal Ab$, hence it 
induces a (contravariant) functor : 
$\text{D}^+(\Cal A)^{\circ}\rightarrow \Cal Ab$ and if one 
applies this functor to the ``long'' complex in 
$\text{D}^+(\Cal A)$ defined in (i) one gets a long exact sequence in 
$\Cal Ab$ (because 
$(X^{\bullet},Y^{\bullet},\text{Con}(u)$$,u,^{\  t}(0,\text{id}_Y),
(\text{id}_{\text{T}X},0))$ is a distinguished triangle in 
$\text{K}(\Cal A)$). 

\vskip 0.1in 

(iv) We also recall that if $I^{\bullet}\in \text{ObK}^+(\Cal A)$ 
consists of injective objects of $\Cal A$ then, for every 
$X^{\bullet}\in \text{ObK}(\Cal A)$, the canonical map 
$\text{Hom}_{\text{K}(\Cal A)}(X^{\bullet},I^{\bullet})
\longrightarrow 
\text{Hom}_{\text{D}(\Cal A)}(X^{\bullet},I^{\bullet})$ is bijective. 

\subhead{3. Example }\endsubhead 
(a) Consider (as in [3] par.3) the short exact sequence in 
$\Lambda \text{-mod}$ : 
$$0\longrightarrow \underline k{\otimes}_kV\longrightarrow   
(\Lambda /({\Lambda}_+)^2)(1)\longrightarrow 
\underline k(1)\longrightarrow 0$$ 
let $w : \underline k(1)\rightarrow 
\text{T}(\underline k{\otimes}_kV)$ be the morphism in 
$\text{D}^b(\Lambda \text{-mod})$ defined in (2)(i) and let 
$\nu=\text{T}^{-1}w : \text{T}^{-1}\underline k(1)\rightarrow 
\underline k{\otimes}_kV$. If one applies L to the short exact 
sequence one gets a $semi$-$split$ short exact sequence in 
$\text{C}(\text{Coh}\Bbb P(V))$. Applying (2)(ii) one derives easily 
that $\text{L}(\nu )$ is the canonical injection : 
${\Cal O}_{\Bbb P}(-1)\rightarrow 
{\Cal O}_{\Bbb P}{\otimes}_kV$ (recall that the module structure 
of $(\Lambda /({\Lambda}_+)^2)(1)$ differs by sign from the module 
structure of $\Lambda /({\Lambda}_+)^2$ ). 

\vskip 0.1in

(b) Dually, consider the short exact sequence in 
$\Lambda \text{-mod}$ : 
$$0\longrightarrow \underline k(-1)\longrightarrow 
(\Lambda /({\Lambda}_+)^2) \spcheck(-1)\longrightarrow 
\underline k{\otimes}_kV^*\longrightarrow 0$$ 
and let $\varepsilon : \underline k{\otimes}_kV^*\rightarrow 
\text{T}\underline k(-1)$ be the morphism in 
$\text{D}^b(\Lambda \text{-mod})$ defined in (2)(i). 
Then $\text{L}(\varepsilon )$ is the canonical epimorphism : 
${\Cal O}_{\Bbb P}{\otimes}_kV^*\rightarrow {\Cal O}_{\Bbb P}(1)$. 

\vskip 0.1in

In the next proposition we gather some well-known properties of the 
category $\Lambda \text{-mod}$, stated in [2]. We include a sketch 
of proof for the reader's convenience.          

\proclaim{4. Proposition } $({\fam0 i})$ If 
$N\in {\fam0 Ob}(\Lambda $-${\fam0 mod})$ and $a\in \Bbb Z$ then the 
map $:$ 
$${\fam0 Hom}_{\Lambda {\fam0 -mod}}(N,\Lambda \spcheck(a))
\longrightarrow 
{\fam0 Hom}_k(N_{-a},\Lambda \spcheck(a)_{-a})=(N_{-a})^*,\   
f\mapsto f_{-a}$$
is bijective. In particular, $\Lambda \spcheck(a)$ is an injective 
object of $\Lambda $-${\fam0 mod}$. 

$({\fam0 ii})$ $\Lambda $-${\fam0 mod}$ has enough injective objects. 

$({\fam0 iii})$ In $\Lambda $-${\fam0 mod}$ $:$ free $\Rightarrow $ 
injective. 

$({\fam0 iv})$ Every $N\in {\fam0 Ob}(\Lambda $-${\fam0 mod})$ has a 
decomposition $:$  
$$N\simeq \Lambda (a_1)\oplus ...\oplus \Lambda (a_m)\oplus N^0$$ 
with $m\in \Bbb N,\  a_1\geq ...\geq a_m$ integers and $N^0$ 
annihilated by ${\fam0 soc}(\Lambda )={\Lambda}_{n+1}$. Moreover, 
$m,a_1,..., a_m$ and $N^0$ $($up to isomorphism$)$ are unique. 

$({\fam0 v})$ In $\Lambda $-${\fam0 mod}$ $:$ projective $\Leftrightarrow $  
free $\Leftrightarrow $ finite direct sum of $\Lambda $-modules of 
the form $\Lambda \spcheck(a)$ $\Leftrightarrow $ injective.
\endproclaim

\demo{Proof } (i) Let 
$f\in \text{Hom}_{\Lambda \text{-mod}}(N,\Lambda \spcheck(a))$. If 
$b>a$, $y\in N_{-b}$ and $\omega \in {\Lambda}_{b-a}$ then 
$f_{-b}(y)\cdot \omega =f_{-a}(y\cdot \omega)$. One can use now the 
fact that the pairing : 
$\Lambda \spcheck(a)_{-b}\times {\Lambda}_{b-a}\rightarrow 
\Lambda \spcheck(a)_{-a}=k$ is perfect. 

(ii) $N$ can be embedded into : 
$\mathop{\oplus}\limits _aN_{-a}{\otimes}_k\Lambda \spcheck(a)$. 

(iii) One can easily show that : 
$\Lambda \simeq \Lambda \spcheck(-n-1)$. 

(iv) For the existence of the decomposition, let $y\in N$ be a 
homogeneous element (let's say, of degree $-a$) not annihilated by 
$\text{soc}(\Lambda )$. Then $y\Lambda \simeq \Lambda (a)$. 
By (ii), $y\Lambda $ is injective in $\Lambda \text{-mod}$ hence it 
is a direct sumand of $N$. One concludes by induction on 
$\text{dim}_kN$. 

For the uniqueness, observe firstly that 
$N\cdot \text{soc}(\Lambda )\simeq \underline k(a_1-n-1)\oplus 
...\oplus \underline k(a_m-n-1)$. This proves the uniqueness of $m$ 
and $a_1,...,a_m$. Assume, now, that one has an isomorphism : 
$$\varphi  : \Lambda (b_1)^{r_1}\oplus ...\oplus 
\Lambda (b_p)^{r_p}\oplus N^0\overset \sim \to\longrightarrow 
\Lambda (b_1)^{r_1}\oplus ...\oplus \Lambda (b_p)^{r_p}\oplus N^1$$ 
with $b_1>...>b_p$ and $N^0,\  N^1$ annihilated by 
$\text{soc}(\Lambda )$. Applying $-\cdot \text{soc}(\Lambda )$ one 
derives that the component of $\varphi $ : 
$\Lambda (b_1)^{r_1}\rightarrow \Lambda (b_1)^{r_1}$ is an 
isomorphism. By a well known trick (about matrices of 2$\times $2=4 
blocks with invertible left upper block) it follows that : 
$$\Lambda (b_2)^{r_2}\oplus ...\oplus \Lambda (b_p)^{r_p}\oplus 
N^0\simeq \Lambda (b_2)^{r_2}\oplus ...\oplus \Lambda (b_p)^{r_p}\oplus 
N^1$$ 
and one concludes by induction on $\text{dim}_kN$. 

(v) Every projective or injective object of $\Lambda \text{-mod}$ is 
a direct sumand of a free object (for injective by the proof of (ii)). 
Now one can apply (iv). 
\qed\enddemo 

\proclaim{5. Lemma } Let $P^{\bullet}\in 
{\fam0 ObC}^-(\Lambda $-${\fam0 mod})$ 
be a complex bounded to the right of free objects of 
$\Lambda $-${\fam0 mod}$. Then the complex ${\fam0 L}(P^{\bullet})$ 
is acyclic. 
\endproclaim 
\demo{Proof } By definition, 
$\text{L}(P^{\bullet})=\text{s}(X^{\bullet \bullet})$ for a double 
complex $X^{\bullet \bullet}$ with 
$X^{p,\bullet}=\text{L}(P^p)$. By (1)(iii), the columns of 
$X^{\bullet \bullet}$ are acyclic bounded complexes. Now, 
$\text{s}(X^{\bullet \bullet})$ is the direct limit of the complexes 
$\text{s}({\sigma}_I^{\geq -p}X^{\bullet \bullet})$, $p\geq 0$, where 
$({\sigma}_I^{\geq -p}X^{\bullet \bullet})^{ij}=X^{ij}$ 
for $i\geq -p$ and $=0$ for 
$i<-p$. ${\sigma}_I^{\geq -p}X^{\bullet \bullet}$ is a 
``first quadrant'' type double complex (i.e., $\exists i_0,\  j_0$ 
such that its $(i,j)$-component is $0$ for $i<i_0$ and, also, for 
$j<j_0$) with acyclic columns, hence 
$\text{s}({\sigma}_I^{\geq -p}X^{\bullet \bullet})$ is acyclic. 
\qed\enddemo 

The next result, which is the key point of this paper, is a 
generalization of the Remark 3 after theorem 2 in [2]. Its proof 
can be easily reduced to the particular case 
$K^{\bullet}=\underline k$ of the remark in [2]. In [2], the remark 
is a consequence of the main result. Here we reverse the order : 
we prove directly the (general version of the) remark and then we 
show that it immediately implies the main result of [2]. 

\proclaim{6. Lemma } Let 
$I^{\bullet}\in {\fam0 ObC}(\Lambda $-${\fam0 mod})$ be an acyclic complex 
of injective $(\Leftrightarrow $ free$)$ objects of 
$\Lambda $-${\fam0 mod}$. For $p\in \Bbb Z$, let 
$Z^p:={\fam0 Ker}d_I^p$. Then $:$  

$({\fam0 a})$ $\  \forall p\in \Bbb Z$, the canonical morphism 
${\fam0 T}^{-p}Z^p\rightarrow I^{\bullet}$ induces a quasi-isomorphism $:$  
${\fam0 L}({\fam0 T}^{-p}Z^p)\rightarrow{\fam0 L}(I^{\bullet})$. 

$({\fam0 b})$ 
$\  \forall K^{\bullet}\in{\fam0 ObC}^b(\Lambda $-${\fam0 mod})$, the 
canonical map $:$ 
$${\fam0 Hom}_{{\fam0 K}(\Lambda)}(K^{\bullet},I^{\bullet})\longrightarrow 
{\fam0 Hom}_{{\fam0 D}(\Bbb P)}({\fam0 L}(K^{\bullet}),
{\fam0 L}(I^{\bullet}))$$ 
is an isomorphism of $k$-vector spaces. 
\endproclaim 

\demo{Proof } (a) Let ${\sigma}^{\geq p}I^{\bullet}$ be the ``stupid'' 
truncation of $I^{\bullet}$ defined by 
$({\sigma}^{\geq p}I^{\bullet})^i=I^i$ for $i\geq p$ and $=0$ for $i<p$. 
The morphism $\text{T}^{-p}Z^p\rightarrow I^{\bullet}$ factorizes as 
$\text{T}^{-p}Z^p\overset \text{qis}\to\longrightarrow 
{\sigma}^{\geq p}I^{\bullet}\rightarrow I^{\bullet}$. One has an exact 
sequence of complexes : 
$$0\longrightarrow{\sigma}^{\geq p}I^{\bullet}\longrightarrow 
I^{\bullet}\longrightarrow{\sigma}^{<p}I^{\bullet}\longrightarrow 0.$$ 

By (5), $\text{L}({\sigma}^{<p}I^{\bullet})$ is acyclic. It follows 
that $\text{L}({\sigma}^{\geq p}I^{\bullet})\rightarrow 
\text{L}(I^{\bullet})$ is a quasi-isomorphism. 

\vskip 0.1in

(b) Let $a:=\text{min}\{ i\in \Bbb Z\mid K^i\neq 0\} $ and 
$b:=\text{min}\{ j\in \Bbb Z\mid K_j^a\neq 0\} $. Then one has a short 
exact sequence : 
$$0\longrightarrow {K'}^{\bullet}\longrightarrow K^{\bullet}\longrightarrow 
{K''}^{\bullet}\longrightarrow 0$$ 
with ${K''}^{\bullet}=\text{T}^{-a}(K_b^a{\otimes}_k\underline k(-b))$. 
Using (2)(iii) and (i) and the Five Lemma one can easily reduce the proof, 
by induction on $\underset i\to\sum \text{dim}_kK^i$, to the case 
$K^{\bullet}=\text{T}^p\underline k(q),\  p,q\in \Bbb Z$, and this case 
reduces immediately to the case $p=q=0$. 

In the case $K^{\bullet}=\underline k$, using (2)(iv) and the fact that 
$\text{T}Z^{-1}\rightarrow{\sigma}^{\geq -1}I^{\bullet}$ and 
$\text{L}(\text{T}Z^{-1})\rightarrow \text{L}(I^{\bullet})$ are 
quasi-isomorphisms one gets isomorphisms : 
$$\text{Hom}_{\text{K}(\Lambda )}(\underline k,I^{\bullet})=
\text{Hom}_{\text{K}(\Lambda )}(\underline k,{\sigma}^{\geq -1}I^{\bullet})
\overset \sim \to\longrightarrow 
\text{Hom}_{\text{D}(\Lambda )}(\underline k,{\sigma}^{\geq -1}I^{\bullet})
\overset \sim \to\longleftarrow $$
$$\overset \sim \to\longleftarrow 
\text{Hom}_{\text{D}(\Lambda )}(\underline k,\text{T}Z^{-1}),$$ 
$$\text{Hom}_{\text{D}(\Bbb P)}(\text{L}(\underline k),
\text{L}(I^{\bullet}))\overset \sim \to\longleftarrow 
\text{Hom}_{\text{D}(\Bbb P)}(\text{L}(\underline k),
\text{L}(\text{T}Z^{-1})).$$ 
It follows that it suffices to prove that the map : 
$$\text{Hom}_{\text{D}(\Lambda )}(\underline k,
\text{T}Z^{-1})\longrightarrow 
\text{Hom}_{\text{D}(\Bbb P)}(\text{L}(\underline k),
\text{L}(\text{T}Z^{-1}))$$ 
is an isomorphism of $k$-vector spaces. We shall prove that, 
$\  \forall N\in \text{Ob}(\Lambda \text{-mod})$: 
$$\text{Hom}_{\text{D}(\Lambda )}(\underline k,
\text{T}^pN)\overset \sim \to\longrightarrow 
\text{Hom}_{\text{D}(\Bbb P)}(\text{L}(\underline k),
\text{L}(\text{T}^pN)),\  \forall p\geq 1. \tag{6.1} $$ 
The proof of (6.1) is based on the following : 
\vskip 0.1in
\bf{Claim:}\  \rm  
$\text{Hom}_{\text{D}(\Lambda )}(\underline k,
\text{T}^p\underline k(a))\overset \sim \to\longrightarrow 
\text{Hom}_{\text{D}(\Bbb P)}(\text{L}(\underline k),
\text{L}(\text{T}^p\underline k(a))),\  \forall p\geq 0,
\  \forall a\in \Bbb Z.$
\vskip 0.1in
Assuming the Claim, for the moment, we prove (6.1) by induction on 
$\text{dim}_kN$. The initial case $\text{dim}_kN=1$ follows from the 
Claim. For the induction step, let 
$a:=\text{min}\{ i\mid N_i\neq 0\} $. One has an exact sequence : 
$0\longrightarrow N'\longrightarrow N\longrightarrow N''\longrightarrow 
0,\  $ 
with $N''=N_a{\otimes}_k\underline k(-a)$. Using the considerations from 
(2)(i), the induction hypothesis for $N'$, the Claim for $N''$ and 
the Five Lemma one gets immediately (6.1). 

Finally, let us prove the Claim. One has : 
$$\text{Hom}_{\text{D}(\Lambda )}(\underline k,
\text{T}^p\underline k(a))\simeq 
\text{Ext}_{\Lambda \text{-mod}}^p(\underline k,\underline k(a)),$$ 
$$\CD \text{Hom}_{\text{D}(\Bbb P)}(\text{L}(\underline k),
\text{L}(\text{T}^p\underline k(a)))=
\text{Hom}_{\text{D}(\Bbb P)}({\Cal O}_{\Bbb P},
\text{T}^{p+a}{\Cal O}_{\Bbb P}(-a))\simeq 
\text{Ext}_{{\Cal O}_{\Bbb P}}^{p+a}({\Cal O}_{\Bbb P},
{\Cal O}_{\Bbb P}(-a)){\simeq}\\
\simeq \text{H}^{p+a}{\Cal O}_{\Bbb P}(-a).\endCD $$ 

$\underline k$ has an injective resolution in $\Lambda \text{-mod}$ : 
$$0\longrightarrow \underline k\longrightarrow 
\Lambda \spcheck\longrightarrow 
V^*{\otimes}_k\Lambda \spcheck(1)\longrightarrow 
...\longrightarrow 
S^iV^*{\otimes}_k\Lambda \spcheck(i)\longrightarrow ...$$
with differential 
$d=\mathop{\sum}\limits _{i=0}^n(X_i\cdot -)\otimes 
(-\cdot e_i)_{\Lambda \spcheck}$. It follows that both sides of the Claim 
are 0 for $p+a\neq 0$ (assuming, of course, $p\geq 0$). It remains to 
show that : 
$$\text{Hom}_{\text{D}(\Lambda )}(\underline k,
\text{T}^p\underline k(-p))\overset \sim \to\longrightarrow 
\text{Hom}_{\text{D}(\Bbb P)}({\Cal O}_{\Bbb P},
{\Cal O}_{\Bbb P}(p)),\  \forall p\geq 0. \tag{6.2} $$ 

Consider the morphism (in $\text{D}(\Lambda )$) 
$\varepsilon : \underline k{\otimes}_kV^*\rightarrow 
\text{T}\underline k(-1)$ from (3)(b). Since 
$\text{L}(\varepsilon )$ is the canonical morphism 
${\Cal O}_{\Bbb P}{\otimes}_kV^*\rightarrow {\Cal O}_{\Bbb P}(1)$, 
$\text{L}(\text{T}^{p-1}\varepsilon (-p+1))$ is the canonical morphism 
${\Cal O}_{\Bbb P}(p-1){\otimes}_kV^*\rightarrow {\Cal O}_{\Bbb P}(p)$. 
Using the commutative diagram : 
$$\CD 
\text{Hom}_{\text{D}(\Lambda )}(\underline k,
\text{T}^{p-1}\underline k(-p+1){\otimes}_kV^*) @>>> 
\text{Hom}_{\text{D}(\Bbb P)}({\Cal O}_{\Bbb P},
{\Cal O}_{\Bbb P}(p-1){\otimes}_kV^*) \\
@VVV @VVV \\
\text{Hom}_{\text{D}(\Lambda )}(\underline k,
\text{T}^p\underline k(-p)) @>>> 
\text{Hom}_{\text{D}(\Bbb P)}({\Cal O}_{\Bbb P},{\Cal O}_{\Bbb P}(p)) 
\endCD $$
one proves easily, by induction on $p\geq 0$, that the morphism in 
(6.2) is surjective, hence it is an isomorphism since both sides are 
isomorphic over $k$ to $S^pV^*$. 
\qed\enddemo 
  
\proclaim{7. Theorem } $({\fam0 Bernstein-Gel'fand-Gel'fand})$ 

$({\fam0 a})$ The functor ${\fam0 L} : \Lambda $-${\fam0 mod}\rightarrow 
{\fam0 D}^b({\fam0 Coh}\Bbb P(V))$ is essentially surjective. 

$({\fam0 b})$ If $N, N'\in {\fam0 Ob}(\Lambda $-${\fam0 mod})$ then 
the map $:$  
$${\fam0 Hom}_{\Lambda {\fam0 -mod}}(N',N)\longrightarrow 
{\fam0 Hom}_{{\fam0 D}^b(\Bbb P)}({\fam0 L}(N'),{\fam0 L}(N))$$
is surjective and its kernel consists of the morphisms factorizing 
through a free $(\Leftrightarrow $ injective $)$ object of 
$\Lambda $-${\fam0 mod}$.
\endproclaim 

\demo{Proof } We firstly prove the second assertion. 
\vskip 0.1in 
(b) Let $0\rightarrow N\rightarrow I^0\rightarrow I^1\rightarrow 
...$ be an injective resolution of $N$ in $\Lambda $-mod and 
$...\rightarrow I^{-2}\rightarrow I^{-1}\rightarrow N\rightarrow 0$ 
a free resolution. Glue them in order to get an acyclic complex 
$I^{\bullet}$ consisting of injective ($\Leftrightarrow $ free) 
objects. By (6)(a), $\text{L}(N)\rightarrow \text{L}(I^{\bullet})$ is 
a quasi-isomorphism and by (6)(b) the bottom horizontal arrow of the 
following commutative diagram : 
$$\CD 
\text{Hom}_{\Lambda \text{-mod}}(N',N) @>>> 
\text{Hom}_{\text{D}^b(\Bbb P)}(\text{L}(N'),\text{L}(N)) \\ 
@VVV @VV\wr V \\ 
\text{Hom}_{\text{K}(\Lambda )}(N',I^{\bullet}) @>\sim >> 
\text{Hom}_{\text{D}(\Bbb P)}(\text{L}(N'),\text{L}(I^{\bullet})) 
\endCD $$ 
is an isomorphism. The left vertical arrow of the diagram is surjective 
and its kernel consists of the morphisms factorizing through 
$I^{-1}$. 

\vskip 0.1in

(a) We observe, firstly, that if 
$K^{\bullet}\in \text{ObC}^b(\Lambda $-mod) then 
$\exists N\in \text{Ob}(\Lambda $-mod) such that 
$\text{L}(K^{\bullet})\simeq \text{L}(N)$ in $\text{D}^b(\Bbb P)$. 
Indeed, consider a quasi-isomorphism 
$u : K^{\bullet}\rightarrow J^{\bullet}$ (resp., 
$v : P^{\bullet}\rightarrow K^{\bullet}$) with 
$J^{\bullet}\in \text{ObC}^+(\Lambda $-mod) (resp., 
$P^{\bullet}\in \text{ObC}^-(\Lambda $-mod)) consisting of injective 
(resp., free) objects. Then $I^{\bullet}:=\text{Con}(u\circ v)$ is an 
acyclic complex consisting of injective ($\Leftrightarrow $ free) 
objects of $\Lambda $-mod. Using the short exact sequence : 
$$0\longrightarrow J^{\bullet}\longrightarrow I^{\bullet}\longrightarrow 
\text{T}P^{\bullet}\longrightarrow 0$$ 
and applying (5) to $\text{T}P^{\bullet}$ one derives that 
$\text{L}(J^{\bullet})\rightarrow \text{L}(I^{\bullet})$ is a 
quasi-isomorphism hence $\text{L}(K^{\bullet})\rightarrow 
\text{L}(I^{\bullet})$ is a quasi-isomorphism. On the other hand, 
by (6)(a), one has a quasi-isomorphism 
$\text{L}(Z^0)\rightarrow \text{L}(I^{\bullet})$. Consequently, 
$\text{L}(K^{\bullet})\simeq \text{L}(Z^0)$ in $\text{D}(\Bbb P)$. 

Let now ${\Cal F}^{\bullet}\in \text{ObC}^b(\text{Coh}\Bbb P(V))$. 
Let $p:=\text{max}\{ i\in \Bbb Z\mid {\Cal F}^i\neq 0\} $ and let 
$u : {\sigma}^{<p}{\Cal F}^{\bullet}\rightarrow 
\text{T}^{-p+1}{\Cal F}^p$ be the morphism defined by 
$d_{\Cal F}^{p-1} : {\Cal F}^{p-1}\rightarrow {\Cal F}^p$. Then 
${\Cal F}^{\bullet}=\text{Con}(u)$. Assume there exist 
$N',N''\in \text{Ob}(\Lambda $-mod) and isomorphisms in 
$\text{D}^b(\Bbb P)$ 
$\psi \  : \text{L}(N'')\overset \sim \to\rightarrow 
{\sigma}^{<p}{\Cal F}^{\bullet}$ , 
$\varphi \  : \text{L}(N')\overset \sim \to\rightarrow 
\text{T}^{-p+1}{\Cal F}^p$. By (b), 
$\exists f\in \text{Hom}_{\Lambda \text{-mod}}(N'',N')$ such that 
$\text{L}(f)={\varphi}^{-1}\circ u\circ {\psi}$. Then 
${\Cal F}^{\bullet}\simeq \text{L}(\text{Con}(f))$ in 
$\text{D}^b(\Bbb P)$, hence, by the above observation, 
$\exists N\in \text{Ob}(\Lambda $-mod) such that 
${\Cal F}^{\bullet}\simeq \text{L}(N)$. 

By induction on the length of ${\Cal F}^{\bullet}$, one can now 
reduce the proof to the case when ${\Cal F}^{\bullet}$ has only 
one non-zero term. By Serre's results from [6], any coherent 
sheaf on $\Bbb P(V)$ has a finite resolution with finite direct 
sums of invertible sheaves ${\Cal O}_{\Bbb P}(a)$. By induction 
on the length of this resolution, one reduces the proof, as above, 
to the case when ${\Cal F}^{\bullet}=\text{T}^p{\Cal O}_{\Bbb P}(a)$. 
But $\text{T}^p{\Cal O}_{\Bbb P}(a)=\text{L}(\text{T}^{p+a}
\underline k(-a))$. 
\qed\enddemo 

\proclaim{8. Corollary } 
$({\fam0 [2]\  Remark\  3\  after\  Theorem\  1})$ 

For every ${\Cal F}^{\bullet}\in {\fam0 ObC}^b({\fam0 Coh}(\Bbb P(V))$ 
there exists $N\in {\fam0 Ob}(\Lambda $-${\fam0 mod})$ annihilated by 
${\fam0 soc}(\Lambda )$ such that 
${\Cal F}^{\bullet}\simeq {\fam0 L}(N)$ in 
${\fam0 D}^b({\fam0 Coh}\Bbb P(V))$. Moreover, $N$ is unique up to 
isomorphism. 
\endproclaim       
 
\demo{Proof } The existence of $N$ follows from (7)(a) and (4)(iv). 
Let $N'$ be another such $\Lambda $-module. By (7)(b), there exists a 
morphism $u : N'\rightarrow N$ in $\Lambda $-mod such that 
$\text{L}(u) : \text{L}(N')\rightarrow \text{L}(N)$ is an 
isomorphism in $\text{D}^b(\Bbb P)$ (i.e., it is a 
quasi-isomorphism). By (7)(b) again, there exists 
$v : N\rightarrow N'$ such that $\text{L}(v)$ is the inverse of 
$\text{L}(u)$ in $\text{D}^b(\Bbb P)$. By the last part of (7)(b), 
there exists a free object $P$ of $\Lambda $-mod such that 
$\text{id}_N-u\circ v$ factorizes as 
$N\overset f\to\rightarrow P\overset g\to\rightarrow N$. 

The submodule of $P$ consisting of the elements annihilated by 
$\text{soc}(\Lambda )$ is 
$P\cdot {\Lambda}_+$, hence 
$f(N)\subseteq P\cdot {\Lambda}_+$, hence 
$\text{Im}(\text{id}_N-u\circ v)\subseteq N\cdot{\Lambda}_+$. 
Using the exterior algebra version of the graded NAK, one derives 
that $u\circ v$ is surjective, hence it is an isomorphism because 
$N$ is a finite dimensional $k$-vector space. Similarly, 
$v\circ u$ is an isomorphism. Consequently, $u$ is an isomorphism. 
\qed\enddemo 

\subhead{9. Definition }\endsubhead 
Let $N\in \text{Ob}(\Lambda $-mod) {\it annihilated by}\rm 
$\  \text{soc}(\Lambda )$. Consider a minimal free resolution of 
$N$ in $\Lambda $-mod : 
$...\rightarrow I^{-2}\rightarrow I^{-1}\rightarrow N\rightarrow 0$. 
Minimality is equivalent to the condition : 
$\text{Im}(I^{-p-1}\rightarrow I^{-p})\subseteq I^{-p}\cdot 
{\Lambda}_+,\  \forall p\geq 1$. Consider also an injective resolution 
of $N$ in $\Lambda $-mod : 
$0\rightarrow N\rightarrow I^0\rightarrow I^1\rightarrow ...$ such 
that $\text{Im}(I^p\rightarrow I^{p+1})\subseteq I^{p+1}\cdot 
{\Lambda}_+,\  \forall p\geq 0$. To get such a resolution, take a 
minimal free resolution of $N\spcheck$ and dualize it. Glueing the 
two resolutions one gets an acyclic complex $I^{\bullet}$ consisting 
of injective ($\Leftrightarrow $ free) objects of $\Lambda $-mod 
such that $\text{Im}d_I^p\subseteq I^{p+1}\cdot 
{\Lambda}_+,\  \forall p\in \Bbb Z$ (for $p=-1$ this follows from 
the fact that $N\cdot \text{soc}(\Lambda )=(0)$). 

Such a complex $I^{\bullet}$ is called a {\it Tate resolution }\rm of 
$N$. 

\proclaim{10. Theorem } 
$({\fam0 Eisenbud-Fl}$\o ${\fam0 ystad-Schreyer})$ 

Let ${\Cal F}^{\bullet}\in {\fam0 ObC}^b({\fam0 Coh}\Bbb P(V))$, let 
$N$ be the unique $($up to isomorphism$)$ object of 
$\Lambda $-${\fam0 mod}$ annihilated by ${\fam0 soc}(\Lambda )$ with 
${\Cal F}^{\bullet}\simeq {\fam0 L}(N)$ in 
${\fam0 D}^b({\fam0 Coh}\Bbb P(V))$ and let $I^{\bullet}$ be a 
Tate resolution of $N$. Then $:$  
\vskip 0.1in 
$({\fam0 a})$ 
$I^p\simeq \underset i\to\bigoplus 
{\Bbb H}^{p-i}{\Cal F}^{\bullet}(i){\otimes}_k\Lambda \spcheck(i),
\  \forall p\in \Bbb Z$ $($where $\Bbb H$ denotes hypercohomology$)$, 

$({\fam0 b})$ $d_I^p : I^p\rightarrow I^{p+1}$ maps 
${\Bbb H}^{p-i}{\Cal F}^{\bullet}(i){\otimes}_k\Lambda \spcheck(i)$ to 
$\underset{j>i}\to\bigoplus{\Bbb H}^{p+1-j}{\Cal F}^{\bullet}(j)
{\otimes}_k\Lambda \spcheck(j)$ and the component $:\   
{\Bbb H}^{p-i}{\Cal F}^{\bullet}(i){\otimes}_k
\Lambda \spcheck(i)\rightarrow 
{\Bbb H}^{p-i}{\Cal F}^{\bullet}(i+1){\otimes}_k
\Lambda \spcheck(i+1)$ 
of $d_I^p$ is defined $($see $({\fam0 4})({\fam0 i}))$ by the 
multiplication map $:\   
{\Bbb H}^{p-i}{\Cal F}^{\bullet}(i){\otimes}_kV^*\rightarrow 
{\Bbb H}^{p-i}{\Cal F}^{\bullet}(i+1)$ 
$($up to sign$)$. 
\endproclaim 

\demo{Proof of $(10)({\fam0 a})$} (according to Remark 3 after 
Theorem 2 in [2]). 

\vskip 0.1in

If $I^p\simeq 
\underset i\to\oplus \Lambda \spcheck(i)^{{\gamma}_{pi}}$ then 
$\text{soc}(I^p)\simeq 
\underset i\to\oplus \underline k(i)^{{\gamma}_{pi}}$. Taking into account 
that $\text{Im}d_I^q\subseteq I^{q+1}\cdot{\Lambda}_+,\  \forall q\in 
\Bbb Z$, one gets that : 
$$\text{soc}(I^p)_{-i}\simeq 
\text{Hom}_{\Lambda \text{-mod}}(\underline k,I^p(-i))\simeq 
\text{Hom}_{\text{K}(\Lambda )}(\underline k,\text{T}^pI^{\bullet}(-i)).$$ 
On the other hand, by (6) : 
$$\text{Hom}_{\text{K}(\Lambda )}(\underline k,
\text{T}^pI^{\bullet}(-i))\simeq 
\text{Hom}_{\text{D}(\Bbb P)}({\Cal O}_{\Bbb P},
\text{T}^{p-i}{\Cal F}^{\bullet}(i))\simeq 
\Bbb E\text{xt}^{p-i}({\Cal O}_{\Bbb P},{\Cal F}^{\bullet}(i))\simeq $$
$\simeq{\Bbb H}^{p-i}{\Cal F}^{\bullet}(i).$ 
\qed\enddemo 
  
For the proof of (10)(b) we need the following addendum to (2)(iii) : 

\subhead{11. Remark }\endsubhead 
Under the assumptions of (2)(iii), let 
$w : Z^{\bullet}\rightarrow \text{T}X^{\bullet}$ be the morphism in 
$\text{D}^+(\Cal A)$ defined in (2)(i). Then :  
$$\text{Hom}_{\text{K}(\Cal A)}(\text{T}^{-1}w,\text{id}_{\text{T}^pI}) : 
\  \text{Hom}_{\text{K}(\Cal A)}(X^{\bullet},\text{T}^pI^{\bullet})
\rightarrow 
\text{Hom}_{\text{K}(\Cal A)}(\text{T}^{-1}Z^{\bullet},
\text{T}^pI^{\bullet})$$ 
equals $(-1)^p{\partial}^p$ where 
${\partial}^p :\  \text{Hom}_{\text{K}(\Cal A)}(X^{\bullet},
\text{T}^pI^{\bullet})\rightarrow 
\text{Hom}_{\text{K}(\Cal A)}(Z^{\bullet},\text{T}^{p+1}I^{\bullet})$ 
is the ``classical'' connecting morphism associated to the short exact 
sequence of complexes of abelian groups : 
$$0\longrightarrow 
\text{Hom}^{\bullet}(Z^{\bullet},I^{\bullet})\longrightarrow 
\text{Hom}^{\bullet}(Y^{\bullet},I^{\bullet})\longrightarrow 
\text{Hom}^{\bullet}(X^{\bullet},I^{\bullet})\longrightarrow 0.$$    

\demo{Proof } ${\partial}^p$ is defined as follows : 
let $f : X^{\bullet}\rightarrow \text{T}^pI^{\bullet}$ be a morphism 
of complexes. Lift every 
$f^i : X^i\rightarrow I^{i+p}$ to a morphism 
$g^i : Y^i\rightarrow I^{i+p}$. Then the morphism of complexes 
$(d_I^{i+p}\circ g^i-(-1)^pg^{i+1}\circ 
d_Y^i)_{i\in \Bbb Z} : Y^{\bullet}\rightarrow 
\text{T}^{p+1}I^{\bullet}$ vanishes on $X^{\bullet}$ hence induces 
a morphism of complexes 
${\partial}^p(f) : Z^{\bullet}\rightarrow \text{T}^{p+1}I^{\bullet}$ 
(in fact, to be rigorous, one has to take homotopy classes). 

We have to prove that the diagram : 
$$\CD 
\text{T}^{-1}\text{Con}(u) @>\pretend{(\text{id}_X,0)}\haswidth
{(-1)^p\text{T}^{-1}{\partial}^p(f)}>> 
X^{\bullet} \\ 
@V(0,\text{T}^{-1}v)VV @VVfV \\ 
\text{T}^{-1}Z^{\bullet} @>(-1)^p\text{T}^{-1}{\partial}^p(f)>> 
\text{T}^pI^{\bullet} 
\endCD $$ 
is homotopically commutative. One can use the homotopy operators 
$h^i:=(0,g^{i-1}) : (\text{T}^{-1}\text{Con}(u))^i=
X^i\oplus Y^{i-1}\rightarrow I^{i+p-1}=(\text{T}^pI^{\bullet})^{i-1}$. 
\qed\enddemo 

\demo{Proof of $(10)({\fam0 b})$} The first assertion follows from the 
fact that $\text{Im}d_I^p\subseteq I^{p+1}\cdot {\Lambda}_+$. For the 
second assertion we consider the morphism (in $\text{D}^b(\Lambda $-mod)) 
$\nu  : \text{T}^{-1}\underline k(1)\rightarrow 
\underline k{\otimes}_kV$ from (3)(a). By (6), the map : 
$$\text{Hom}_{\text{K}(\Lambda )}(\nu ,\text{id}) :\    
\text{Hom}_{\text{K}(\Lambda )}(\underline k{\otimes}_kV,
\text{T}^pI^{\bullet}(-i))\longrightarrow 
\text{Hom}_{\text{K}(\Lambda )}(\text{T}^{-1}\underline k(1),
\text{T}^pI^{\bullet}(-i))$$ 
can be identified to the map : 
$$\text{Hom}_{\text{D}(\Bbb P)}(\text{L}(\nu ),\text{id}) :\  
\text{Hom}_{\text{D}(\Bbb P)}({\Cal O}_{\Bbb P}{\otimes}_kV,
\text{T}^{p-i}{\Cal F}^{\bullet}(i))\rightarrow 
\text{Hom}_{\text{D}(\Bbb P)}({\Cal O}_{\Bbb P}(-1),
\text{T}^{p-i}{\Cal F}^{\bullet}(i))$$ 
and this one can be identified to the multiplication map : 
${\Bbb H}^{p-i}{\Cal F}^{\bullet}(i){\otimes}_kV^*\rightarrow 
{\Bbb H}^{p-i}{\Cal F}^{\bullet}(i+1)$. 

We want now to explicitate $\text{Hom}_{\text{K}(\Lambda )}(\nu ,
\text{id})$. Let $\xi \in {\Bbb H}^{p-i}{\Cal F}^{\bullet}(i)$, 
$\lambda \in V^*$ and let 
$f :\  \underline k{\otimes}_kV\rightarrow I^p(-i)$ be the morphism 
defined by 
$\xi \otimes \lambda  :\  \underline k{\otimes}_kV\rightarrow 
{\Bbb H}^{p-i}{\Cal F}^{\bullet}(i){\otimes}_k(\Lambda \spcheck)_0$. 
Then $f$ can be lifted to the morphism  
$g : (\Lambda /({\Lambda}_+)^2)(1)\rightarrow I^p(-i)$ sending 
$\hat 1\in (\Lambda /({\Lambda}_+)^2)(1)_{-1}$ to 
$-\xi \otimes \lambda \in{\Bbb H}^{p-i}{\Cal F}^{\bullet}(i)
{\otimes}_k(\Lambda \spcheck)_{-1}$. 

Using (11) (and the explicit description of ${\partial}^p$ from its 
proof) one derives that 
$\text{Hom}_{\text{K}(\Lambda )}(\nu ,\text{id})$ can be identified to : 
$$(-1)^{p-1}d_I^p \mid \  
{\Bbb H}^{p-i}{\Cal F}^{\bullet}(i){\otimes}_k
\Lambda \spcheck(i)_{-i-1}\longrightarrow 
{\Bbb H}^{p-i}{\Cal F}^{\bullet}(i+1){\otimes}_k  
\Lambda \spcheck(i+1)_{-i-1}.\qed $$ 
\enddemo 

One can easily deduce from (10) the Lemma of Castelnuovo-Mumford. More 
important, Eisenbud et al.[3] show that (10) implies the results of 
A.A. Beilinson [1]. We close the paper by briefly explaining this, in 
terms of the present approach. 

\proclaim{12. Theorem } $({\fam0 Beilinson})\  $ Let 
${\Cal F}^{\bullet}\in{\fam0 ObC}^b({\fam0 Coh}\Bbb P(V))$. Then 
${\Cal F}^{\bullet}$ is isomorphic in 
${\fam0 D}^b({\fam0 Coh}\Bbb P(V))$ to a complex $C^{\bullet}$ with 
$C^p=\underset i\to\bigoplus{\Bbb H}^{p+i}{\Cal F}^{\bullet}(-i)
{\otimes}_k{\Omega}_{\Bbb P}^i(i),\  \forall p\in \Bbb Z$ and also to 
a complex ${C'}^{\bullet}$ with 
${C'}^p=\underset i\to\bigoplus{\Cal O}_{\Bbb P}(-i){\otimes}_k
{\Bbb H}^{p+i}({\Cal F}^{\bullet}\otimes{\Omega}_{\Bbb P}^i(i)),\  
\forall p\in \Bbb Z$. 
\endproclaim 

\demo{Proof } (according to [3] (6.1) and (8.11)). 
\vskip 0.1in 
Let $N$ and $I^{\bullet}$ be as in the statement of (10). Recall, from 
(6)(a), that $\text{L}(N)\rightarrow \text{L}(I^{\bullet})$ is a 
quasi-isomorphism. By definition, 
$\text{L}(I^{\bullet})=\text{s}(X^{\bullet \bullet})$ for a certain 
double complex $X^{\bullet \bullet}$ with 
$X^{pq}={\Cal O}_{\Bbb P}(q){\otimes}_kI_q^p$. 
\vskip 0.1in 
(I) In order to prove the first assertion, one takes 
$C^{\bullet}:=\text{Ker}(X^{\bullet ,0}\rightarrow X^{\bullet ,1})$. 
Taking into account that 
$\text{Ker}(\text{L}(\Lambda \spcheck(-i))^0\rightarrow 
\text{L}(\Lambda \spcheck(-i))^1)={\Omega}_{\Bbb P}^i(i)$, the formula 
for $C^p$ follows from (10)(a). 

It remains to show that 
$C^{\bullet}\rightarrow \text{s}(X^{\bullet \bullet})$ is a 
quasi-isomorphism. It decomposes as 
$C^{\bullet}\rightarrow 
\text{s}({\sigma}_{II}^{\geq 0}X^{\bullet \bullet})\rightarrow 
\text{s}(X^{\bullet \bullet})$, where 
$({\sigma}_{II}^{\geq 0}X^{\bullet \bullet})^{ij}:=X^{ij}$ for 
$j\geq 0$ and $=0$ for $j<0$. Since the columns of $X^{\bullet \bullet}$ 
are acyclic, 
$C^{\bullet}\rightarrow \text{s}({\sigma}_{II}^{\geq 0}X^{\bullet \bullet})$ 
is a quasi-isomorphism. On the other hand, one has a short exact 
sequence : 
$$0\longrightarrow \text{s}({\sigma}_{II}^{\geq 0}X^{\bullet \bullet})
\longrightarrow \text{s}(X^{\bullet \bullet})\longrightarrow 
\text{s}({\sigma}_{II}^{<0}X^{\bullet \bullet})\longrightarrow 0$$ 
hence it suffices to prove that 
$\text{s}({\sigma}_{II}^{<0}X^{\bullet \bullet})$ is acyclic. 

For $p\in \Bbb Z$, let 
$a(p):=\text{min}\{ q\in \Bbb Z\mid I_q^p\neq 0\} $. Since 
$I^{p-1}\rightarrow I^p\rightarrow Z^{p+1}\rightarrow 0$ is a 
minimal free presentation, it follows that $a(p-1)>a(p),\  
\forall p\in \Bbb Z$. One deduces that {\it the rows of }\rm 
$X^{\bullet \bullet}$ {\it are bounded to the left}\rm . 

Now, $\text{s}({\sigma}_{II}^{<0}X^{\bullet \bullet})$ is the direct 
limit of the complexes 
$\text{s}({\sigma}_{II}^{\geq -p}{\sigma}_{II}^{<0}X^{\bullet \bullet})$ 
for $p\geq 1$. But 
${\sigma}_{II}^{\geq -p}{\sigma}_{II}^{<0}X^{\bullet \bullet}$ is a 
``first quadrant'' type double complex with acyclic rows hence its 
associated simple complex is acyclic. 
\vskip 0.1in 
(II) Let us prove the second assertion. One takes the subcomplex 
$J^{\bullet}$ of $I^{\bullet}$ defined by 
$J^p:=\underset{i\geq 0}\to\bigoplus 
{\Bbb H}^{p-i}{\Cal F}^{\bullet}(i){\otimes}_k\Lambda \spcheck(i)$. 
One has $J^p=0$ for $p<<0$ and $J^p=I^p$ for $p>>0$ hence, by (5), 
$\text{L}(J^{\bullet})\rightarrow \text{L}(I^{\bullet})$ is a 
quasi-isomorphism. 

Now, $\text{L}(J^{\bullet})=\text{s}(Y^{\bullet \bullet})$ for a certain 
double complex $Y^{\bullet \bullet}$ with splitting rows 
$Y^{\bullet ,q},\  q\in \Bbb Z$ and with columns 
$Y^{p,\bullet}=\text{L}(J^p)=0$, for $p<<0$. According to a general 
lemma about such double complexes (see [3] (3.5)), 
$\  \text{s}(Y^{\bullet \bullet})$ is homotopically equivalent to a 
complex ${C'}^{\bullet}$ whose ``linear part'' is 
$\text{L}(\text{H}^{\bullet}(J^{\bullet}))$ (where 
$\text{H}^{\bullet}(J^{\bullet})$ is the complex with $p$th term 
$\text{H}^p(J^{\bullet})$ and with all the differentials equal to 0). 
In particular : 
${C'}^p\simeq \underset i\to\bigoplus 
{\Cal O}_{\Bbb P}(-i){\otimes}_k\text{H}^{p+i}(J^{\bullet})_{-i}$. 
But $J_{-i}^{\bullet}=0$ for $i<0$ and 
$J_{-i}^{\bullet}=I_{-i}^{\bullet}$ for $i>n$ hence 
(since $I^{\bullet}$ is acyclic) 
$\text{H}^q(J^{\bullet})_{-i}=0$ for $i<0$ and for $i>n$,
$\  \forall q\in \Bbb Z$. On the other hand, for 
$0\leq i\leq n$ and $q\in \Bbb Z$ : 
$J_{-i}^q\simeq 
\text{Hom}_{\Lambda \text{-mod}}((\Lambda /({\Lambda}_+)^{i+1})(i),I^q)$ 
hence $\text{H}^q(J^{\bullet})_{-i}\simeq 
\text{Hom}_{\text{K}(\Lambda )}((\Lambda /({\Lambda}_+)^{i+1})(i),
\text{T}^qI^{\bullet})$. One can now apply (6), taking into account 
that $\text{L}((\Lambda /({\Lambda}_+)^{i+1})(i))\simeq 
({\Omega}_{\Bbb P}^i(i))^*$ in 
$\text{D}^b(\Bbb P)$. 
\qed\enddemo

\enddocument